\input xy
\xyoption{all}

\normalbaselineskip=1.6\normalbaselineskip\normalbaselines
\magnification=1200
\outer\def\beginsection#1\par{\vskip0pt plus0.2\vsize\penalty-250
	\vskip0pt plus-.2\vsize\bigskip\vskip\parskip
	\message{#1}\leftline{\bf#1}\nobreak\smallskip\noindent}
\def\ra{\rightarrow}
\def\mZ{{\bf Z}}
\def\mQ{{\bf Q}}
\def\mC{{\bf C}}
\def\mF{{\bf F}}
\def\Spec{{\rm Spec}}
\def \bs{\bigskip}
\def\Gal{\mathop{\rm Gal}\nolimits}

\def\mod{{\rm mod}}

\def\Stab{\mathop{\rm Stab}\nolimits}
\def\cl{\mathop{\rm cl}\nolimits}
\def\Ker{\rm Ker}
\def\mtr#1{{\overline{#1}}}
\def\unr{{\rm unr}}
\def\Tor{\mathop{\rm Tor}\nolimits}
\def\ref#1{{\it #1}}
\def\Zar{\mathop{\rm Zar}\nolimits}

\def\Id{\mathop{\rm Id}\nolimits}
\def\pr{'}

\centerline{\bf A note on the Manin-Mumford conjecture}
\medskip
\centerline{\bf Damian Roessler
\footnote{$^\dagger$}{\rm 
CNRS, Institut de math\'ematiques de Jussieu, 
Universit\'e Paris 7, Case Postale 7012, 2, place Jussieu, 
75251 PARIS CEDEX 05, 
FRANCE} 
}

{\bf Abstract.} 
In [PR1], R. Pink and the author gave a short 
proof of the Manin-Mumford conjecture, which was inspired 
by an earlier model-theoretic proof by Hrushovski. 
The proof given in [PR1] uses a difficult unpublished 
ramification-theoretic 
result of Serre. It is the purpose of this note 
to show how the proof given in [PR1] can be modified 
so as to circumvent the reference to Serre's result. 
J. Oesterl\'e and R. Pink contributed several simplifications 
and shortcuts to this note. 
\bs

\beginsection 0. Introduction.

Let $A$ be an abelian variety defined over an algebraically closed 
field $L$ of characteristic $0$ and let 
$X$ be a closed subvariety. If $G$ is an abelian group, write 
$\Tor(G)$ for the group of elements of $G$ which are of finite 
order. A closed subvariety of $A$ whose irreducible 
components are translates of abelian subvarieties of $A$
 by torsion points will be called a torsion subvariety.
The Manin-Mumford conjecture is the following statement:

{\it The Zariski closure of $\Tor(A(L))\cap X$ is 
a torsion subvariety.}

\noindent
This was first proved by Raynaud in [R]. 
In [PR1], R. Pink and the author gave a new proof of this 
statement, which was inspired by an earlier model-theoretic 
proof given by Hrushovski in [H]. The interest of this proof 
is the fact that it relies almost entirely on classical algebraic geometry 
and is quite short. Its only non elementary input is a  
ramification-theoretic result of 
Serre. The proof of this result is not published and relies 
(see [Se] (pp.~33--34, 56--59)) on deep theorems of 
Faltings, Nori and Raynaud.  
In this note, we show how the reference to Serre's result in [PR1] 
can be replaced by a reference 
to a classical result in the theory of formal groups (see Th. 4 (a)). 

The structure of the paper is as follows. 
For the convenience of the reader, the text has been written 
so as to be logically independent of [PR1]. In particular, no 
knowledge of [PR1] is necessary to read it.  
Section 1 recalls various classical results on abelian varieties 
and also contains two less well-known, but elementary propositions 
(Prop. 1 and Prop. 3) whose proofs can be found elsewhere but 
for which we have included short proofs to make the text 
more self-contained.  The reader 
is encouraged to proceed directly to section 2, which contains 
a complete proof of the Manin-Mumford conjecture and 
to refer to the results listed in section  1 
as needed. 

{\bf Notations.} w.r.o.g. is a shortening of {\it without restriction 
of generality}; if $X$ is closed subvariety of an abelian 
variety $A$ defined over an algebraically closed field $L$ 
of characteristic $0$, then we write $\Stab(X)$ for the 
stabiliser of $X$; this is a 
closed subgroup of $A$ such that $\Stab(X)(L):=\{a\in A(L)|a+X=X\}$;  
it has the same field of definition as $X$ and $A$; if $p$ is a prime 
number and $G$ is an abelian group, we write $\Tor^p(G)$ for 
the set of elements of $\Tor(G)$ whose order is prime to $p$ and 
$\Tor_p(G)$ for the set of elements of $\Tor(G)$ whose order is a power 
of $p$. 

{\bf Acknowledgments.} We want to thank J. Oesterl\'e 
for his interest and for suggesting some simplifications 
in the proofs of [PR1] (see [Oes]) which have inspired some 
of the proofs given here. Also, the proof 
of Prop. 3 in its present form is due to him 
(see the explanations before the proof).
I am also very grateful to R. Pink, who carefully read several 
versions of the text and suggested many 
improvements and simplifications. 
In particular, Prop. 6 was suggested by him.  

\beginsection 1. Preliminaries. 

\proclaim Lemma 0. 
Let $L\subseteq L\pr$ be algebraically closed fields of 
characteristic $0$. 
Let $A$ be an abelian variety defined over 
$L$ and let $X$ be a closed $L$-subvariety of $A$. Then:
\item{(a)} $X$ is a torsion subvariety of $A$ iff $X_{L\pr}$ is 
a torsion subvariety of $A_{L'}$;
\item{(b)} the Manin-Mumford 
conjecture holds for $X$ in $A$ iff it holds for $X_{L\pr}$ 
in $A_{L\pr}$. 

\noindent{\it Proof:} we first prove (a). To prove 
the equivalence of the two conditions, we only need to prove 
the sufficiency of the second one. The latter is 
a consequence of the fact that the morphism $\pi:A_{L\pr}\ra A$ is 
faithfully flat and that any torsion point and any
abelian subvariety of $A_{L\pr}$ has a model in $A$ 
(see [Mi] (Cor. 20.4, p. 146)). To prove (b),  
let $Z:=\Zar(\Tor(A(L))\cap X)$ (resp. 
$Z\pr:=\Zar(\Tor(A(L\pr))\cap X_{L\pr})$). Using again the fact that 
 any torsion point in $A_{L\pr}$ has a model in $A$ and 
that $\pi$ is faithfully flat, we see that 
$\pi^{-1}(\Tor(A(L))\cap X)=\Tor(A(L\pr))\cap X_{L\pr}$. 
From this and the fact that the morphism $\pi$ is open 
([EGA] (IV, 2.4.10)), we get a set-theoretic equality 
$\pi^{-1}(Z)=Z\pr$. Since $\pi$ is radicial, the underlying 
set of $\pi^*(Z):=Z_{L\pr}$ is $\pi^{-1}(Z)$ 
([EGA] (I, 3.5.10)). Since $Z_{L\pr}$ is 
reduced ([EGA] (IV, 4.6.1)), we thus have an 
equality of closed subschemes $Z_{L\pr}=Z\pr$. 
Now, by (a), the closed subscheme 
$Z_{L\pr}$ is a torsion subvariety of $A_{L\pr}$ iff  
$Z$ is a torsion subvariety of $A$. 
$\bullet$

\proclaim Proposition 1 (Pink-Roessler). 
Let $A$ be an abelian variety over $\mC$ and 
let $F:A\ra A$ be an isogeny. 
Suppose that the absolute 
value of all the eigenvalues of the pull-back 
map $F^*$ on the first singular cohomology group $H^1(A(\mC),\mC)$ is larger than $1$. 
Then any closed subvariety $Z$ of $A$ such 
that $F(Z)=Z$ is a torsion subvariety.

\noindent
The following proof can be found in [PR1] (Remark after Lemma 2.6).

\noindent{\it Proof:} 
w.r.o.g., we may replace 
$F$ by one of its powers and thus suppose that each 
irreducible component of $Z$ is stable under $F$. 
We may thus suppose that $Z$ is irreducible.  
Notice that $F(\Stab(Z))\subseteq\Stab(Z)$. 
Let us first suppose that $\Stab(Z)=0$. 

Write $\cl(Z)$ for the cycle class of 
$Z$ in $H^*(A(\mC),\mC)$. 
We list the following facts:
\item{(1)} the degree of $F$ 
 is the determinant of the restriction of 
$F^*$ to $H^1(A(\mC),\mC)$;
\item{(2)} each eigenvalue of $F^*$ on 
$H^i(A(\mC),\mC)$ is the product 
of $i$ distinct zeroes (with multiplicities) 
of the characteristic polynomial of $F^*$ on $H^1(A(\mC),\mC)$; 
\noindent Facts (1) and (2) follow from the fact that 
for all $i\geq 0$ there is a natural isomorphism 
$\Lambda^i(H^1(A(\mC),\mC))\simeq H^i(A(\mC),\mC)$ (see [Mu] (p.3, Eq. (4))). 

Now notice that since $\Stab(Z)=0$, the varieties 
$Z+a$, where $a\in\Ker(F)(\mC)$, are pairwise 
distinct. These varieties are thus the irreducible components of 
$F^{-1}(Z)$. Now we compute
$$
\cl(F^*(Z))=\sum_{a\in\Ker(F)}\cl(Z+a)=
\#\Ker(F)(\mC)\cdot\cl(Z)=\deg(F)\cl(Z)
$$
and thus $\cl(Z)$ belongs to the eigenspace of 
the eigenvalue $\deg(F)$ in $H^*(A(\mC),\mC)$. Facts (1), (2) and 
the hypothesis on the eigenvalues 
imply that 
$\cl(Z)\in H^{2\dim(A)}(A(\mC),\mC)$, which in turn 
implies that $Z$ is a point. This point is a torsion 
point since it lies in the kernel of $F-\Id$, which is 
an isogeny by construction.

If $\Stab(Z)\ne 0$, then replace 
$A$ by $A/\Stab(Z)$ and $Z$ by $Z/\Stab(Z)$. 
The isogeny $F$ then induces an isogeny 
on $A/\Stab(Z)$, which stabilises $Z/\Stab(Z)$. 
We deduce that $Z/\Stab(Z)$ is 
a torsion point. This implies that $Z$ is a translate 
of $\Stab(Z)$ by a torsion point and concludes 
the proof.
$\bullet$

\proclaim Corollary 2. 
Let $A$ be an abelian variety over an algebraically closed field $K$ of characteristic $0$.
Let $n\geq 1$ and let $M$ be an $n\times n$-matrix with integer 
coefficients. Suppose that the absolute 
value of all the eigenvalues of $M$ is larger than $1$. 
Then any closed subvariety $Z$ of $A^n$ such that 
 $M(Z)=Z$ is a torsion subvariety.

\noindent{\it Proof:}   
Because of Lemma 0 (a), we may assume 
w.r.o.g. that $K$ is the algebraic closure of 
a field which is finitely generated as a field over $\mQ$. We may thus 
also assume that $K\subseteq\mC$. Prop. 1 then implies 
the result for $Z_\mC$ in $A^n_\mC$ and using Lemma 0 (a) again we 
can conclude. $\bullet$

\proclaim Proposition 3 (Boxall). 
Let $A$ be an abelian variety over a field $K$ of characteristic $0$. 
Let $p>2$ be a prime number and 
let $L:=K(A[p])$ be the extension of $K$ 
generated by the $p$-torsion points of $A$. 
Let $P\in \Tor_p(A(\mtr{K}))$ and suppose that $P\not\in A(L)$.
Then there exists $\sigma\in\Gal(\mtr{L}|L)$ 
such that $\sigma(P)-P\in A[p]\setminus\{0\}$.

\noindent A proof of a variant of Prop. 3 can be found in [B]. 
For the convenience of the reader, we 
reproduce a proof, which is a simplification by 
Oesterl\'e (private communication) of a proof 
due to Coleman and Voloch (see [Vo]).

\noindent{\it Proof:} 
let $n\geq 1$ be the smallest natural number so that $p^nP \in A(L)$. 
For all $i\in\{1,\dots,n\}$, let $P_i = p^{n-i}P$. Let also $\sigma_1$ be an element of 
$\Gal(\mtr{L}|L)$ such that $\sigma_1(p^{n-1}P) \ne p^{n-1}P$.
Furthermore, let $\sigma_i:=\sigma_1^{p^{i-1}}$ and $Q_i:=\sigma_i(P_i) - P_i$. 

First, notice that we have $pQ_1=\sigma_1(p^nP)-p^nP=0$ and $Q_1=
\sigma_1(p^{n-1}P)-p^{n-1}P\not=0$, hence $Q_1\in A[p]\backslash\{0\}$.
We shall prove by
induction on $i\geq 1$ that $Q_i=Q_1$ if $i\leq n$. 
This will prove the proposition, since $Q_n=\sigma_n(P)-P$. 

So assume that $Q_i=Q_1$ for some $i<n$. We have
$p^2(\sigma_i-1)(P_{i+1})=p(\sigma_i-1)(P_i)=pQ_i=0$. Since any
$p$-torsion point of $A$ is fixed by $\sigma$, and hence by $\sigma_i$,
we also have $p (\sigma_i-1)^2(P_{i+1})=0$ and
$(\sigma_i-1)^3(P_{i+1})=0$.  The binomial formula shows that, in the
ring of polynomials ${\bf Z}[T]$,  $T^p$ is congruent to $1+p(T-1)$
modulo  the ideal generated by $p(T-1)^2$ and $(T-1)^3$ (notice that
$p\not=2$ !). We thus have
$(\sigma_i^p-1)(P_{i+1})=p (\sigma_i-1)(P_{i+1})= (\sigma_i-1)(P_i)$,
id est  $Q_{i+1}=Q_i$. This completes the induction 
on $i$.
$\bullet$

Suppose now that $K$ is a finite extension of 
$\mQ_p$, for some prime number $p$ and let $K^\unr$ be its maximal 
unramified extension. Let $k$ be the 
residue field of $K$.  
Suppose that $A$ is an abelian variety over $K$ which has good reduction at 
the unique non-archimedean place of $K$. Denote by $A_0$ 
the corresponding special fiber, which is an abelian 
variety over $k$.

\proclaim Theorem 4. 
\item{(a)} The kernel of the homomorphism
$$
\Tor(A(K^\unr))\ra A_0(\mtr{k})
$$
induced by the reduction map is a finite $p$-group.
\item{(b)} The equality $\Tor^p(A(K^\unr))=\Tor^p(A(\mtr{K}))$ holds.

\noindent{\it Proof:} for statement (b), see [Mi] (Cor. 20.8, p. 147). 
Statement (a), which is more difficult to prove, follows from general 
properties of formal groups over $K$. See [Oes2] (Prop. 2.3 (a)) 
for the proof.   
$\bullet$

Let now $\phi\in\Gal(\mtr{k}|k)$ be the arithmetic 
Frobenius map. 

\proclaim Theorem 5 (Weil).
There is a monic polynomial $Q(T)\in\mZ[T]$ 
with the following properties:
\item{(a)} $Q(\phi)(P)=0$ for all $P\in A_0(\mtr{k})$;
\item{(b)} the complex roots of $Q$ have absolute value 
$\sqrt{\#k}$.

\noindent{\it Proof:} see [We].$\bullet$

\beginsection 2. Proof of the Manin-Mumford conjecture.

\proclaim Proposition 6. 
Let $A$ be an abelian variety over a field $K_0$ 
that is finitely generated as a field over $\bf Q$. Then for almost all 
prime numbers $p$, there exists an embedding of $K_0$ 
into a finite extension $K$ of ${\bf Q}_p$, such that $A_K$ has good 
reduction at the unique non-archimedean place of $K$.

\noindent{\it Proof:} 
since by assumption $K_0$ has finite transcendence degree 
over $\mQ$, there is a finite map 
$$\Spec\ K_0\ra \Spec\ \mQ(X_1,\dots,X_d),$$ for 
some $d\geq 0$ (notice that $d=0$ is allowed).  
Let $V\ra {\bf A}^d_\mZ$ be the normalisation of the affine space ${\bf A}^d_\mZ$ 
in $K_0$. The scheme $V$ is integral, normal 
and has $K_0$ as a field of rational functions. Furthermore, 
$V$ is finite and surjective onto ${\bf A}^d_\mZ$. 
There is an open subset $B\subseteq V$ and an abelian 
scheme ${\cal A}\ra B$, whose generic fiber is $A$. Choose 
$B$ sufficiently small so that its image $f(B)$ is open 
and so that $f^{-1}(f(B))=B$ (this can be achieved by replacing 
$B$ by $f^{-1}({\bf A}^d_\mZ\backslash f(V\backslash B))$). 
Let $U:=f(B)$. This accounts for the square on the left of the diagram (*) 
below.

Now notice that $U(\mQ)\not=\emptyset$, since 
${\bf A}^d(\mQ)$ is dense in ${\bf A}^d_{\mQ}$ and 
$U\cap {\bf A}^d_{\mQ}$ is open and not empty. 
Thus, for almost all prime numbers $p$, we have $U(\mF_p)\not=\emptyset$. 
Let $p$ be a prime number with this property. 
Let $P\in U(\mF_p)$ and let $a_1,\dots,a_d\in\mF_p$ be 
its coordinates. Choose 
as well elements $x_{1},\dots, x_{d}\in{{\mQ_p}}$ 
which are algebraically independent over $\mQ$. The elements 
$x_1,\dots,x_d$ remain algebraically 
independent if we  
replace some $x_i$ by ${1\over x_i}$ so we may suppose that 
$\{x_1,\dots,x_d\}\subseteq{\cal O}_{\mQ_p}$. Notice also that 
any element of the residue field $\mF_p$ of ${\cal O}_{\mQ_p}$ 
is the reduction $\mod\ p$ of an element of $\mZ\subseteq{\cal O}_{\mQ_p}$. 
Furthermore, the elements $x_1,\dots,x_d$ remain algebraically independent 
if some $x_i$ is replaced by $x_i+m$, where $m$ is 
an integer. Hence, we may also suppose that $x_i\ \mod\ p=a_i$ for 
all $i\in\{1,\dots,d\}$. 
The choice of the $x_i$ induces a morphism $e:\Spec\ {\cal O}_{\mQ_p}\ra 
{\bf A}^d_\mZ$, which by construction sends the generic point of 
$\Spec\ {\cal O}_{\mQ_p}$ on the generic point of 
${\bf A}^d_\mZ$ and hence of $U$ and sends the special point of $\Spec\ {\cal O}_{\mQ_p}$ 
on $P\in U(\mF_p)$. Hence $e^{-1}(U)=\Spec\ {\cal O}_{\mQ_p}$. 
This accounts for the lowest square in (*).

The middle square in (*) is obtained by taking the fibre product of 
$U\ra B$ and $\Spec\ {\cal O}_{\mQ_p}\ra B$.  The morphism $B_1\ra\Spec\ {\cal O}_{\mQ_p}$ is 
then also finite and surjective.

To define the arrows in the triangle next to it, consider 
a reduced irreducible component 
$B_1\pr$ of $B_1$ which dominates $\Spec\ {\cal O}_{\mQ_p}$. 
This exists, because the morphism $B_1\ra\Spec\ {\cal O}_{\mQ_p}$ is 
dominant. The morphism $B_1\pr\ra\Spec\ {\cal O}_{\mQ_p}$ will then also 
be finite and 
will thus correspond to a finite (and hence integral) extension of integral 
rings. Let $K$ be the function field 
of $B_1\pr$, which is a finite extension of 
$\mQ_p$; the ring associated to $B_1\pr$ is by construction 
included in the integral closure ${\cal O}_K$ of 
${\cal O}_{\mQ_p}$ in $K$ and the arrow 
$\Spec\ {\cal O}_K\ -->\ B_1$ is defined 
by composing the morphism induced by this inclusion with 
the closed immersion $B_1\pr\ra B_1$.

The morphism $\Spec\ K\ra\Spec\ {\mQ_p}$ has been
 implicitly defined in 
the last paragraph and the morphisms 
$\Spec\ {\mQ}_p\ra\Spec\ {\cal O}_{\mQ_p}$ and 
$\Spec\ K\ra\Spec\ {\cal O}_{K}$ are the obvious ones. 

We have a commutative diagram (*):

\xymatrix
{
\Spec\ K_0\ar@{}[dr]|{\rm Cart.}\ar@{=>}[d]\ar[r] & B\ar@{=>}[d]\ar@{}[dr]|{\rm Cart.} & 
B_1\ar[l]\ar@{=>}[d] & \ar@{-->}[l]\Spec\ {\cal O}_{K}\ar@{=>}[dl] &
\Spec\ K\ar[l]\ar@{=>}[d]\\
\Spec\ \mQ(X_1,\dots, X_d)\ar[r] & U\ar[d]\ar@{}[dr]|{\rm Cart.} & 
\ar[l]\Spec\ {\cal O}_{\mQ_p}\ar@{=}[d] & & \Spec\ \mQ_p\ar[ll]\\
& {\bf A}^d_\mZ & \Spec\ {\cal O}_{\mQ_p}\ar[l] & &  
}

{\it The single-barreled continuous arrows ($\ra$) represent dominant maps; 
the double-barreled continuous ones ($\Rightarrow$) represent finite and dominant maps; 
all the schemes in the diagram apart from $B_1$ are integral;  
the cartesian squares carry the label "Cart.".}

Now notice that the map 
$\Spec\ K\ra B$ obtained by composing the connecting morphisms 
sends $\Spec\ K$ on 
the generic point of $B$; to see this notice that 
the maps $\Spec\ K\ra\Spec\ {\cal O}_K$, 
$\Spec\ {\cal O}_K\Rightarrow\Spec\ {\cal O}_{\mQ_p}$ and 
$\Spec\ {\cal O}_{\mQ_p}\ra U$ are all dominant; hence 
$\Spec\ K$ is sent on 
the generic point of $U$; since $B\ra U$ is a finite 
map, this implies that $\Spec\ K$ is sent on the generic point of 
$B$. 

Thus the map $\Spec\ K\ra B$ induces a field 
extension $K|K_0$.  
Furthermore, as we have seen, $K$ is a finite extension of 
$\mQ_p$ and by construction, the abelian variety 
$A_K$  
is  
the generic fiber of the abelian scheme  
${\cal A}\times_B \Spec\ {\cal O}_K$. In other 
words ${A}_K$ is an abelian variety defined over $K$ which 
has good reduction at the unique non-archimedean place of 
$K$.$\bullet$

Next, we shall consider the following situation. 
Let $p>2$ be a prime number and let 
$K$ be a finite extension of $\mQ_p$. Let 
$k$ be its residue field. Let $A$ be an abelian 
variety over $K$.  
Suppose that $A$ has good reduction
at the unique non-archimedean place of $K$. Let 
 $A_0$ be the corresponding special fiber, which is an abelian 
variety over $k$. 

Recall that $K^\unr$ refers to the maximal unramified extension of $K$. 
Let $\phi\in\Gal(\mtr{k}|k)$ be the arithmetic Frobenius map 
and let $\tau\in\Gal(K^\unr|K)$ be its canonical lift. 

\proclaim Proposition 7. 
Let $X$ be a closed $K$-subvariety of $A$.  
Then the Zariski closure of $X_{\mtr{K}}\cap\Tor(A(K^\unr))$ is 
a torsion subvariety.

\noindent{\it Proof:} w.r.o.g. 
we may suppose that $\Tor(A(K^\unr))$ is dense in 
$X_{\mtr{K}}$ (otherwise, replace $X$ by the
 natural model of $\Zar(X_{\mtr{K}}\cap\Tor(A(K^\unr)))$ over $K$). 
By Th. 4 (a), the kernel of the reduction 
homomorphism $\Tor(A(K^\unr))\ra A_0(\mtr{k})$ is a finite $p$-group. 
Let $p^r$ be its cardinality and let $Y:=p^r\cdot X$. 
Let
$Q(T):=T^n-(a_{n}T^{n-1}+\dots+a_0)\in\mZ[T]$ be the polynomial 
provided by Th. 5 (i.e. the characteristic polynomial of 
$\phi$ on $A_0(\mtr{k})$). Let $F$ be the matrix
$$
\pmatrix {0&1&\ldots&0&0\cr
\vdots&\vdots&&\vdots&\vdots\cr
0&0&\ldots&0&1\cr
a_0&a_1&\ldots&a_{n-2}&a_{n-1}\cr}
$$
\smallskip
\noindent
For any $a\in A(K^\unr)$, write 
$u(x):=(x,\tau(x),\tau^2(x),\dots,\tau^{n-1}(x))\in A^n(K^\unr)$.
Let $\widetilde{Y}:=\Zar(\{u(a) | a\in (p^r\cdot\Tor(A(K^\unr)))\cap Y_{\mtr{K}}\})$. 
Th. 5 (a) and Th. 4 (a) imply that
$$F(u(a))=u(\tau(a))$$
 for all 
$a\in p^r\cdot\Tor(A(K^\unr))$. Furthermore, 
by construction, 
$$\tau(p^r\cdot\Tor(A(K^\unr)))\subseteq 
p^r\cdot\Tor(A(K^\unr)).$$
 Hence $F(\widetilde{Y})=\widetilde{Y}$. 
Now Th. 5 (b) implies that the absolute value of the 
eigenvalues of the matrix $F$ are larger than $1$ 
and Cor. 2 then implies that $\widetilde{Y}$ is a torsion subvariety of 
$A_{\mtr{K}}$.
The variety $Y_{\mtr{K}}$ is the projection of $\widetilde{Y}$ on the first 
factor and is thus also a torsion subvariety. Finally, 
this implies that $X_{\mtr{K}}$ is a torsion subvariety. 
$\bullet$

\proclaim Proposition 8.  
Let $X$ be a closed $K$-subvariety of $A$. 
Then the Zariski closure of $X_{\mtr{K}}\cap\Tor(A(\mtr{K}))$ is 
a torsion subvariety.

\noindent{\it Proof:} we may suppose w.r.o.g. that $K=K(A[p])$, 
that $X$ is 
geometrically irreducible 
and that $X_{\mtr{K}}\cap\Tor(A(\mtr{K}))$ is dense in $X_{\mtr{K}}$. 
We shall first suppose that $\Stab(X)=0$. 
Let $x\in X_{\mtr{K}}\cap\Tor(A(\mtr{K}))$ 
and suppose that $x\not\in A(K^\unr)$. 
Write $x=x^p+x_p$, where $x^p\in \Tor^p(A(\mtr{K}))$ and 
$x_p\in\Tor_p(A(\mtr{K}))$. By Th. 4 (b) 
$x^p\in A(K^\unr)$ and thus $x_p\not\in A(K^\unr)$.  
By Prop. 3, there exists $\sigma\in 
\Gal(\mtr{K}|K^\unr)$ such that 
$$
\sigma(x_p)-x_p=\sigma(x)-x\in A[p]\setminus\{0\}.
$$
Now notice that for all $y\in X(\mtr{K})$ and all 
$\tau\in \Gal(\mtr{K}|K^\unr)$, we have $\tau(y)\in  X(\mtr{K})$. 
Hence if the set $\{x\in X_{\mtr{K}}\cap\Tor(A(\mtr{K}))|x\not\in 
A(K^\unr)\}$ is dense in $X_{\mtr{K}}$ then 
$\Stab(X)(\mtr{K})$ contains a element of $A[p]\setminus\{0\}$. 
Since $\Stab(X)=0$, we deduce that 
the set $\{x\in X_{\mtr{K}}\cap\Tor(A(\mtr{K}))|x\not\in 
A(K^\unr)\}$  is not dense in $X_{\mtr{K}}$ and thus 
the set $X_{\mtr{K}}\cap\Tor(A(K^\unr))$ is dense in $X_{\mtr{K}}$. 
Prop. 7 then implies that $X_{\mtr{K}}$ is a torsion point. 
If $\Stab(X)\ne 0$, then we may apply the same reasoning 
to $X/\Stab(X)$ and $A/\Stab(A)$ to conclude that 
$X_{\mtr{K}}$ is a translate of $\Stab(X)_{\mtr{K}}$ by a torsion point. 
$\bullet$

We shall now prove the Manin-Mumford conjecture. 
Let the terminology of the introduction hold. 
By Lemma 0 (b), we may assume w.r.o.g. that 
$L$ is the algebraic closure of a field $K_0$ that is 
finitely generated as a field over $\bf Q$ and that 
$A$ (resp. $X$) has a model $\bf A$ (resp. $\bf X$) 
over $K_0$. By Prop. 6, 
there is 
an embedding of $K_0$ into a field $K$, with 
the following properties: $K$ is a finite 
extension of $\mQ_p$, where $p$ is a prime number larger than $2$ 
and ${\bf A}_K$ has good reduction at 
the unique non-archimedean place of $K$. 
Prop. 8 now implies that the Manin-Mumford conjecture holds 
for ${\bf X}_{\mtr{K}}$ in 
${\bf A}_{\mtr{K}}$ and using Lemma 0 (b) we 
deduce that it holds for $X$ in $A$.  

\noindent
{\bf Remark.} Let the notation of 
the introduction hold. Prop. 3. {\it alone} implies the statement 
of the Manin-Mumford conjecture, with $\Tor(A(L))$ replaced 
by $\Tor_p(A(L))$, for any prime number $p>2$. To see this, 
we may w.r.o.g. assume that $X$ is irreducible and 
that $\Tor_p(A(L))\cap X$ is dense in $X$. By 
an easy variant of Lemma 0 (b), we may w.r.o.g. assume 
that $L$ is the algebraic closure of a field $K$ that is 
finitely generated as a field over $\bf Q$ and that 
$A$ (resp. $X$) has a model $\bf A$ (resp. $\bf X$) 
over $K$. Finally, we may assume w.r.o.g. that 
$K=K({\bf A}[p])$. 
Suppose first that $\Stab(X)=0$. 
By the same argument as above, the set 
$\{a\in \Tor_p(A(L))|a\not\in {\bf A}(K),\ a\in 
X\}$ 
is not 
dense in $X$. Hence the set 
$\{a\in \Tor_p(A(L))|a\in {\bf A}(K),\ a\in X\}$ must be 
dense in $X$; the theorem of Mordell-Weil  
(for instance) implies that this set is finite and thus $X$ 
consists of a single torsion point. If $\Stab(X)\ne 0$, then 
we deduce by the same reasoning that $X/\Stab(X)$ is a torsion 
point in $A/\Stab(X)$ and hence $X$ is a translate of $\Stab(X)$ 
by a torsion point. 
This proof of a special case of the Manin-Mumford conjecture is 
outlined in [B] (Remarque 3, p. 75).

\bigskip
\centerline{\bf References.}
\bigskip
\noindent
[B] Boxall, J. Sous-vari\'et\'es alg\'ebriques de vari\'et\'es
semi-ab\'eliennes sur un corps fini in {\it Number Theory, Paris 1992-3},
S. David, ed., London Math. Soc. lecture note series {\bf 215}, 69--89,
Cambridge Univ. Press, 1995.
\medskip
\noindent
[EGA] Grothendieck, A. \'El\'ements de g\'eom\'etrie alg\'ebrique. 
 {\it Inst. Hautes \'Etudes Sci. Publ. Math.} {\bf 4, 8, 11, 17, 20, 
24, 28, 32} (1960-1967).
\medskip
\noindent
[H] Hrushovski, E. The Manin-Mumford conjecture and the
model theory of difference fields. {\it Ann. Pure Appl. Logic}
{\bf 112} (2001), no. 1, 43--115.
\medskip
\noindent
[Mi] Milne, J. S. Abelian varieties.  {\it Arithmetic geometry (Storrs, Conn., 1984)},  103--150, Springer, New York, 1986.
\medskip
\noindent
[Mu] Mumford, D. Abelian varieties. {\it Tata Institute of Fundamental Research Studies in Mathematics, No. 5}, Oxford University Press, London, 1970.
\medskip
\noindent
[Oes] Oesterl\'e, J. Lettre \`a l'auteur (20/12/2002).
\medskip
\noindent
[Oes2] Oesterl\'e, J. Courbes sur une vari\'et\'e ab\'elienne 
(d'apr\`es M. Raynaud). S\'eminaire Bourbaki, Vol. 1983/84.  
{\it Ast\'erisque}  No. {\bf 121-122} (1985), 213--224.
\medskip
\noindent
[PR1]  Pink, R., Roessler, D. 
On Hrushovski's proof of the Manin-Mumford conjecture. 
{\it Proceedings of the International Congress of Mathematicians}, Vol. I (Beijing, 2002), 539--546, Higher Ed. Press, Beijing, 2002. 
\medskip
\noindent
[R] Raynaud, M. Sous-vari\'et\'es d'une vari\'et\'e ab\'elienne et
 points de torsion. 
{\it Arithmetic and geometry}, Vol. I, 327--352, Progr. 
Math. {35}, Birkh\"auser Boston, Boston, MA, 1983.
\medskip
\noindent
[Se] Serre, J.-P. Oeuvres, vol. IV (1985-1998). Springer 2000.
\medskip
\noindent
[Vo] Voloch, J.-F. 
Integrality of torsion points on abelian varieties over $p$-adic fields.
{\it Math. Res. Lett.} {\bf 3} (1996), no. 6, 787--791.
\medskip
\noindent
[We] Weil, A. 
Vari\'et\'es ab\'eliennes et courbes alg\'ebriques.
{Hermann} 1948.
\end